\documentclass[letterpaper, 10 pt, conference]{ieeeconf}  

\IEEEoverridecommandlockouts                              
\usepackage{times} 
\usepackage{amsmath,amsfonts} 
\usepackage{amssymb}  
\usepackage{tabularx}
\usepackage{ulem}
\usepackage{dsfont}
\usepackage{graphicx}
\usepackage[dvipsnames,usenames]{color}
\usepackage[colorlinks,%
linkcolor=BrickRed,%
filecolor=BrickRed,%
citecolor=RoyalPurple,%
]{hyperref}
\usepackage{color}
\usepackage{todonotes}

\def\Re{\mathbb{R}}
\def\R{\mathbb{R}}

\def\argmin{\mathop{\text{\rm arg\,min}}}

\def\notes#1{\marginpar{\tiny #1}\typeout{Notes!
Notes!
Notes!
}}
\renewcommand{\notes}[1]{\typeout{notes!}}
\def\FRAC#1#2#3{\genfrac{}{}{}{#1}{#2}{#3}}
\def\half{{\mathchoice{\FRAC{1}{1}{2}}%
{\FRAC{2}{1}{2}}%
{\FRAC{3}{1}{2}}%
{\FRAC{4}{1}{2}}}}

\def\Re{\field{R}}
\def\k{{\sf K}}

\def\eqdef{\mathrel{:=}}

\def\Expect{{\mathbb E}}

\def\R{\mathbb{R}}

\def\eqdef{\mathrel{:=}}

\newtheorem{theorem}{Theorem}

\newtheorem{remark}{Remark}
\newtheorem{proposition}{Proposition}

\def\beq{\begin{eqnarray}} 
\def\bc{\begin{center}} 
\def\be{\begin{enumerate}}
\def\bi{\begin{itemize}} 
\def\bs{\begin{small}}
\def\bS{\begin{slide}}
\def\ec{\end{center}} 
\def\ee{\end{enumerate}}
\def\ei{\end{itemize}}
\def\es{\end{small}}
\def\eS{\end{slide}}
\def\eeq{\end{eqnarray}}

\newcommand{\trace}{\text{Tr}}

\newcommand{\ud}{\,\mathrm{d}}

\def\Re{\mathbb{R}}

\def\argmin{\mathop{\text{\rm arg\,min}}}

\renewcommand{\Re}{\mathbb{R}}

\def\eqdef{\mathbin{:=}}

\def\FRAC#1#2#3{\genfrac{}{}{}{#1}{#2}{#3}}

\newcounter{rmnum}
\newenvironment{romannum}{\begin{list}{{\upshape (\roman{rmnum})}}{\usecounter{rmnum}
\setlength{\leftmargin}{14pt}
\setlength{\rightmargin}{8pt}
\setlength{\itemsep}{2pt}
\setlength{\itemindent}{-1pt}
}}{\end{list}}

\newcounter{anum}

\title{\LARGE \bf
	An Optimal Transport Formulation of Bayes' Law \\for Nonlinear Filtering Algorithms
}

\author{Amirhossein Taghvaei and Bamdad Hosseini
	\thanks{A. Taghvaei is with the Department of Aeronautics \& Astronautics, University of Washington, Seattle {\tt\small amirtag@uw.edu}} 
	\thanks{B. Hoesseini is with the Department of Applied Mathematics, University of Washington, Seattle {\tt\small bamdadh@uw.edu}}%
}

\begin{document}

	\maketitle
	\thispagestyle{empty}
	\pagestyle{empty}

	\begin{abstract}
		This paper presents a variational representation of the Bayes' law using optimal transportation theory.  The variational representation is in terms of the optimal transportation between the joint distribution of the (state, observation) and their independent coupling. By imposing certain structure on the transport map, the solution to the variational problem is used to construct a Brenier-type map that transports the prior distribution to the posterior distribution for any value of the observation signal. 
		The new formulation is used  to derive the optimal transport form of the Ensemble Kalman filter (EnKF) for the discrete-time filtering problem and propose a novel extension of EnKF to the non-Gaussian setting utilizing input convex neural networks. Finally, the proposed methodology is used to derive  the optimal transport form of the feedback particle filler (FPF) in the continuous-time limit, which constitutes  its first variational construction without explicitly using the nonlinear filtering equation or Bayes' law.
	\end{abstract}
	\section{Introduction}
	
		Nonlinear filtering is the problem of computing the conditional distribution of the state of a stochastic dynamical system given historical noisy observations. The critical step in any nonlinear filtering algorithm is the implementation of Bayes' law in order to update the conditional distribution of the 
		state as the new observations arrive.
		Bayes' law gives the conditional distribution of the state $X$ given observations $Y$ according to
	\begin{equation}\label{eq:Bayes}
		P_{X|Y} (x|y)= \frac{P_X(x) P_{Y|X}(y|x)}{P_Y(y)} 
	\end{equation}
	where $P_X$ is the prior distribution assumed on $X$, $P_{Y|X}$ is the likelihood distribution of observations conditioned
		on the state $X$, and $P_Y(y) = \int P_{Y|X}(y|x)P_X(x)\ud x$ is the marginal distribution of $Y$.
		Nonlinear filtering is also accompanied with an update step according to the dynamics of the system. This step is straightforward using
the dynamic model directly. Therefore, the focus of this paper is on the Bayesian update step. 
		
		Numerical implementations of Bayes' law often require a discretization (finite-dimensional approximation)
		of the prior and posterior distributions. An exact computation is not possible except in a few special cases, such as the class of Gaussian distributions or in the finite-state space setting. 
		This motivates Monte-Carlo or particle-based approaches where 
		the posterior distribution is approximated
		with an empirical distribution of samples. To this end, the main task of a particle-based algorithm is to transform a 
	set of particles that represent samples from the prior distribution to particles that represent  the posterior distribution. Simply put, the problem can be stated as:
	\begin{align*}
		\text{given:}&~\{X^1_0,\ldots,X^N_0\} \sim P_X\\
		\text{generate:}&~\{X^1_1,\ldots,X^N_1\} \sim P_{X|Y}
	\end{align*}
	
	A vanilla importance sampling and resampling particle filter carries out this task by first forming a weighted empirical distribution and then resampling from the weighted distribution~\cite{doucet09}: 
	\begin{align*}
		X^i_1 \sim \sum_{j=1}^N w_j \delta_{X^j_0},\quad w_j = \frac{P_{Y|X}(y|X^j_0)}{\sum_{k=1}^N P_{Y|X}(y|X^k_0)}.
	\end{align*}
	Although computationally efficient, this approach performs poorly in high-dimensional problems  due to weight degeneracy
	issues~\cite{rebeschini2015can,beskos2014error,bengtsson08}. 
	In particular, for a  Gaussian prior and  identity observation model, the mean-squared error 
	scales with $\frac{C^d}{N}$ where $d$ is the dimension of the state and $C$ is a positive constant~\cite{taghvaei2020optimal}. 
	This means that in order to keep the same error, the number of particles $N$ should scale
	exponentially with the dimension $d$.
	
	These issues motivated recent efforts in the nonlinear filtering literature to develop numerical algorithms
	based on a controlled system of interacting particles to approximate the posterior
	distribution~\cite{taoyang_TAC12, yang2016, crisan10, reich11, reich2015probabilistic, bergemann2012ensemble,
		daum10,daum2017generalized}. A prominent idea is to view the problem of transforming samples from the prior to
	the posterior from the lens of optimal transportation theory~
	\cite{reich2013nonparametric,reich2019data,taghvaei2020optimal,reich13,AmirACC2016,taghvaei2021optimal}, which has
	also become popular in the Bayesian inference literature~\cite{el2012bayesian,marzouk2016introduction,
		mesa2019distributed,heng2015gibbs,spantini2019coupling,kovachki2020conditional,siahkoohi2021preconditioned}.
	{
		Broadly speaking, the aim of the above methods is to find a  {\it transport map} (be it stochastic or deterministic)
		that transforms the prior distribution to the posterior distribution while minimizing a certain cost.

		The present paper builds on the aforementioned transport-based  works and is closely related to the optimal transport
		formulation of the Ensemble Kalman filter (EnKF)  and feedback particle filter
		(FPF)~\cite{taghvaei2020optimal,taghvaei2021optimal} with a slight, but crucial,  difference in the formulation:
		Prior works are based on optimal transportation from the prior $P_X$ to the posterior $P_{X|Y}$,  which is undesirable
		as the posterior distribution is not available and depends on the observation $Y$.}  Instead,  the new formulation
	involves the optimal transportation between the independent coupling $P_X \otimes P_Y$ and the joint distribution
	$P_{XY}$ which is readily available in a filtering problem. By imposing a block triangular
	structure $(x,y) \mapsto S(x,y)=(T(x,y),y)$ on the transport map $S$, the  component  $T(x,y)$ automatically
	serves as a transport map between the prior and posterior distributions for any realization of the  observation $Y=y$. Such triangular representations of the transport map have also appeared in~\cite{el2012bayesian,spantini2019coupling,kovachki2020conditional,ray2022efficacy,shi2022conditional}. With the  map $T$ at hand, one can generate samples from the posterior distribution $P_{X|Y=y}$ according to
			\begin{equation*}
			X^i_1 = T(X^i_0,y),\quad \text{for}\quad i=1,\ldots,N.
			\end{equation*}
		
	
	{
		This slight difference in the formulation of the transport problem has several implications
		which constitute the contributions of the paper:
		\begin{romannum}
			\item A variational formulation of Bayes' law is introduced that is based on the Kantorovich dual formulation of
			a  optimal transportation problem between $P_X \otimes P_Y$ and $P_{XY}$.
			The solution to this problem is a Brenier-type  map (the map $T$ above) that characterizes the
			posterior.    
			\item The variational formulation is cast as a stochastic optimization problem that requires pairs of samples $\{(X^i,Y^i)\}_{i=1}^N$
			from the joint distribution $P_{XY}$. The pair of samples can be generated using a simulator/oracle  for the observation model,  without the need for analytical model of the likelihood function.
			\item The stochastic optimization problem involves a search over a certain class of convex functions. With a quadratic restriction of the  search domain, the optimal transport formulation of EnKF  is obtained, while a  restriction to the class of input convex neural networks (ICNN)~\cite{AmoXuKol17} leads to novel nonlinear filtering algorithms that generalize EnKF to non-Gaussian settings.  
			\item In the continuous-time limit, the solution to the variational problem is used to recover the optimal transport version of the FPF algorithm~\cite{taghvaei2021optimal}, which constitutes the first variational construction of FPF that does not explicitly rely on the nonlinear filtering equations.     
		\end{romannum}

		The rest of the paper is organized as follows: Necessary background on optimal transportation is
		summarized in Section~\ref{sec:OT}. The variational  formulation of Bayes' law
		is presented in Section~\ref{sec:OT-Bayes}. Computational algorithms are discussed in~\ref{sec:computational}, and the connection to the FPF algorithm is discussed in~\ref{sec:FPF}.} Concluding remarks are given in Section~\ref{sec:conclusion}.

	
	\subsection{Background on optimal transportation theory}\label{sec:OT}
	Given two probability measures  $\mu,\nu$ on $\Re^n$, and a measurable map $T:\Re^n \to \Re^n$, we say $T$ transports $\mu$ to $\nu$, if $T_\sharp \mu = \nu$, where $\sharp$ is the push-forward operator. In probabilistic terms, if $Z$ is a random variable with probability
	law $\mu$, then $T(Z)$ has probability law  $\nu$.
	The set of all transport maps from $\mu$ to $\nu$ is denoted by $\mathcal T(\mu,\nu)$. 
	
	
		The Monge optimal transportation problem with quadratic cost is to select the transport map from $\mu$ to $\nu$ with the least quadratic cost:
		\begin{equation*}
			\min_{T \in \mathcal T(\mu,\nu)}\Expect_{Z \sim \mu} \left[\frac{1}{2}\|T(Z)-Z\|^2 \right].
		\end{equation*}  
		If $\mu$ is absolutely continuous with respect to the Lebesgue measure then the Monge problem has a
		dual formulation, referred to as the Monge-Kantorovich (MK) dual problem \cite[Thm.~2.9]{villani2003}
		\begin{equation*}
			\min_{f \in \text{CVX}(\mu)} \Expect_{Z \sim \mu}[f(Z)] + \Expect_{V \sim \nu}[f^*(V)]
		\end{equation*}    
		where $f^*$ denotes the convex conjugate of  $f$, i.e. $f^*(v) = \sup_{z \in \Re^n} z^Tv - f(z)$, and  $\text{CVX}(\mu)$  denotes the set of all convex and $\mu$-integrable functions on $\mathbb R^n$. Then the celebrated result of Brenier \cite{brenier1987decomposition}
		states that the MK dual problem above has a unique minimizer $\bar{f}$ and that $\bar{T} = \nabla \bar{f}$
		is the solution to the Monge problem. The map $\bar{T}$ is often referred to as the Brenier transport map (solution).

	\section{A variational formulation of Bayes' law}\label{sec:OT-Bayes}
		This section proposes a variational formulation of Bayes' law~\eqref{eq:Bayes},
		characterizing the posterior distribution as the pushforward of the prior via 
  a parametric map.
	We assume that the hidden state $X \in \mathcal \Re^n$ and the observation $Y \in \mathcal \Re^m$.
	In the prior works~\cite{taghvaei2020optimal,taghvaei2021optimal}, the authors considered the problem of
	finding the optimal transport map $x \mapsto T(x)$ that transports the prior distribution $P_X$
	to the posterior distribution $P_{X|Y}$:
	\begin{equation}\label{eq:OT-prior-work}
		\min_{T \in \mathcal T(P_X,P_{X|Y})}~\Expect_{X \sim P_X}[\|T(X)-X \|^2 ].
	\end{equation}
	
	This formulation is not useful for constructing a filtering algorithm as it involves the unknown conditional distribution $P_{X|Y}$ explicitly. {Furthermore, the optimal map $\bar{T}$ has to be re-computed for each observation of $Y$.}
	The key idea that resolves this issue is to instead
	consider a transport problem on the product space $\Re^n \times \Re^m$.
	More precisely, 
	find a map $(x,y) \mapsto S(x,y)=(T(x,y),y)$ that
	transports the independent coupling $P_X \otimes P_Y$ to the joint distribution $P_{XY}$. 
	The structure of the map $S(x,y)$  implies that its first component $T(x,y)$ serves as
	a transport map from $P_X$ to $P_{X|Y}$ yielding the new optimal transportation problem
	{
		\begin{equation}\label{eq:OT-S}
			\min_{ S \in \mathcal T(P_X \otimes P_Y, P_{XY})} \Expect_{(X,Y)\sim P_X \otimes P_Y}[\| T(X,Y)-X\|^2],
		\end{equation}
		where the map $S$ is assumed to have the above parameterization in terms of $T$ and we used the identity
		$\|S(X,Y) - (X,Y)\|^2 = \|T(X,Y)-X\|^2$.  In contrast to the prior formulation~\eqref{eq:OT-prior-work},
		the new formulation contains the posterior distribution implicitly through the joint distribution $P_{XY}$
		which is readily available for the filtering task and without knowledge of a likelihood distribution.}
	
	{
		The optimal transportation problem~\eqref{eq:OT-S} is numerically infeasible to solve as it requires
		search over the set $\mathcal T(P_X \otimes P_Y, P_{XY})$. Instead, following~\cite{peyre2019computational},
		we consider the MK dual problem:
		\begin{equation}\label{eq:OT-f}
			\begin{aligned}
				\min_{f(\cdot,y) \in \text{CVX}(P_X)}  & \mathbb E_{(X,Y) \sim P_X \otimes P_Y }[f( X,Y)]  \\ & \hspace{5ex}+ \mathbb E_{(X,Y) \sim P_{XY} }[f^*(X,Y)], 
			\end{aligned}
		\end{equation}
	where the constraint $f(\cdot;y) \in \text{CVX}(P_X)$ means that $x \mapsto f(x;y)$ is convex  and in $L^1(P_X)$
	for any $y \in \Re^m$. Similarly,  $f^*(x;y)=\sup_{z} z^Tx - f(z;y)$ is the convex conjugate  of $f(\cdot;y)$  for fixed $y$.}
%

The following  theorem, which is a direct consequence of \cite[Thm.~2.3]{carlier2016vector} and can be viewed as a conditional analogue of Brenier's result,
connects the minimizers of \eqref{eq:OT-S} and \eqref{eq:OT-f}. 
\medskip

\begin{theorem}
	Assume $\mathbb E[\|X\|^2]<\infty$ and $P_X$ is absolutely continuous with
		respect to the Lebesgue measure on $\R^n$.
	Then, the MK dual problem \eqref{eq:OT-f} admits an optimal solution $\bar f$ that satisfies
	\begin{equation}
		P_{X|Y=y} =  \nabla_x \bar{f}(\cdot,y)\sharp P_X,\quad \forall y \in \R^m.
	\end{equation}
	Moreover, $\bar{T}(x, y) = \nabla_x \bar{f}(x,y)$ is the solution to \eqref{eq:OT-S}. 
	%
\end{theorem}
\medskip

\begin{proof}(sketch)
The main idea behind the proof is the disintegration of the objective function in~\eqref{eq:OT-f} according to
\begin{align*}
    \int \left[\int f(x,y) \ud P_{X}(x)  + \int f^*(x,y) \ud P_{X|Y}(x|y) \right]\ud P_Y(y), 
\end{align*}
and noting that the term inside brackets is the MK dual problem for transporting from $P_X$ to $P_{X|Y=y}$ for each $y \in \mathbb R^m$. Then, by the application of the Brenier's result to the MK dual problem for each $y$, a minimizer $g_y$ exists such that $P_{X|Y=y} = \nabla {g_y} \sharp P_X$. The proof follows by showing that the function $f$, defined according to $f(x,y):=g_y(x)$ for all $x \in \mathbb R^n$ and $y\in \mathbb R^m$, is the optimal solution for~\eqref{eq:OT-f}. 
\end{proof}

\begin{remark}
	The proposed variational formulation for the posterior distribution is fundamentally different than the existing variational approaches based on Kullback-Leibler (KL) divergence that solve
	\begin{align*}
	& P_{X|Y=y} = \argmin_{Q}\,D(Q\|P_{X|Y=y} ),\\
	&=\argmin_{Q}\, \left\{D(Q\|P_X) - \int \log(P_{Y|X}(y|x))dQ(x)\right\},
	\end{align*}     
	where $D(Q\|P) = \int \log(\frac{dQ}{dP})dQ$ is the KL divergence and Bayes' law was used to obtain
	the second display. This formulation is reminiscent of the JKO time stepping procedure~\cite{jordan1998variational} and has been used extensively in variational  construction of filtering algorithms~\cite{laugesen15,halder2017gradient,halder2018gradient, halderproximal}. However, theoretically, it involves the  Bayes' law explicitly, and computationally, it is not meaningful when the prior distribution is an empirical distribution formed by particles. 
\end{remark}	

\section{Computational algorithms} \label{sec:computational}
The proposed variational formulation is computationally very appealing and gives rise to
novel families of algorithms. In particular, let 
	\begin{equation*}
	J(f) \eqdef \mathbb E_{(X,Y) \sim P_X \otimes P_Y }[f( X,Y)] + \mathbb E_{(X,Y) \sim P_{XY} }[f^*(X,Y)], 
	\end{equation*}
	denote the objective function in~\eqref{eq:OT-f}. The value of the objective function is easily  approximated empirically. In particular, 
	given the ensemble of particles $\{X^i_0\}_{i=1}^N$ that form samples from the prior distribution $P_X$, one generates observations $Y^i_0 \sim P_{Y|X}(\cdot|X^i_0)$ for $i=1,\ldots,N$ so that $\{(X^i_0,Y^i_0)\}_{i=1}^N$ form independent samples from the joint distribution $P_{XY}$. The samples are then used to define the empirical cost
\begin{equation*}
	J^{(N)}(f) \eqdef \frac{1}{N^2}\sum_{i,j=1}^N f(X^i_0,Y^j_0) +\frac{1}{N}\sum_{i=1}^N  f^*(X^i_0,Y^i_0).
\end{equation*}
It is straightforward  to see that $J^{(N)}(f)$ is an unbiased estimator of $J(f)$. Remarkably, it is not necessary to use the analytical  form of the observation likelihood in order to form the estimate. Only a simulator/oracle to generate observations is required.  

The empirical approximation is then utilized to formulate  optimization problems of the form:
\begin{equation}\label{eq:OT-fN}
	\min_{f \in \mathcal F}~J^{(N)}(f) 
\end{equation}
where  $\mathcal F$ is a subset of  functions $f(x;y):\Re^{n} \times \Re^m \to \Re$ such that
$x \mapsto f(x;y)$ is convex in $x$. The solution of the optimization problem, denoted by $f^{(N)}$, forms a numerical approximation of the solution $\bar{f}$ to~\eqref{eq:OT-f} which is used to transport the particles
\begin{equation*}
X^i_1= \nabla_x f^{(N)}(X^i_0,y)
\end{equation*}
for the received realization of the observation $Y=y$. 
 
We discuss restriction of $\mathcal F$ to two class of functions, namely quadratic, and neural networks,
in Section~\ref{sec:EnKF} and~\ref{sec:ICNN} respectively. The latter
is closely related to the Monotone GANs algorithm of \cite{kovachki2020conditional}
and the cWGAN algorithm of \cite{ray2022efficacy}.

{
\subsection{Optimal transport EnKF}\label{sec:EnKF}
Consider the class of quadratic functions in $x$,
\begin{align}
	\mathcal F_{Q} = \Big\{(x;y) \mapsto \frac{1}{2}x^TAx +& x^T( K y + b)\,\Big|\,A \in S^{n}_{+},\nonumber \\&\quad
	K \in \Re^{n \times m}, b \in \Re^n \Big\}\label{eq:FQ}
\end{align}
where $S^n_+$ denotes the set of positive-definite matrices. 
Such functions give rise to linear transport  maps  of the form
\begin{equation}
	\nabla_x f(x,y) = Ax + Ky + b,
\end{equation}
which are sufficient to represent exact transport maps from priors  to posteriors whenever $P_{XY}$ is Gaussian.
Note that in this setting both the prior and the posterior are Gaussian.

Considering problem \eqref{eq:OT-f} for functions  $ f \in \mathcal F_Q$, which in turn is analogous to
optimizing over the parameters 
$\theta:=\{A,K,\tilde{b}\} \in \Theta:=S^n_+ \times \Re^{n\times m} \times \Re^n$
of a quadratic function, yields the optimization problem}
\begin{align}\label{eq:J-FQ}
&\min_{\theta \in \Theta}\,\frac{1}{2}\trace(A\Sigma_{x}) +\frac{1}{2}\trace(A^{-1}\Sigma_{x}) +\frac{1}{2}\trace(A^{-1}K\Sigma_{y}K^T) \nonumber
\\
&\quad  - \trace(A^{-1}\Sigma_{xy}K^T)+\frac{1}{2} (\tilde{b}-m_x)^TA^{-1}(\tilde{b}-m_x)
\end{align}
where $m_x = \mathbb E [X]$, $m_y = \mathbb E [Y]$, $\tilde{b} = b - Am_x -Km_y$,   $\Sigma_x = \mathbb E[(X-m_x)(X-m_x)^T]$, $\Sigma_{y} = \mathbb E[(Y-m_y)(Y-m_y^T)]$  and $\Sigma_{xy} = \Expect[(X-m_x)(Y-m_y)^T]$. 
The following proposition  characterizes the exact solution to \eqref{eq:J-FQ}.  
\begin{proposition}\label{prop:quadratic}
Assume $\Sigma_X$ 
is  positive definite. Then, the optimization problem~\eqref{eq:J-FQ} is convex and admits the unique solution 
\begin{subequations}
	\begin{align}
		\tilde{b} &=  m_x \\\label{eq:optimal-K}
		K &= \Sigma_{xy}\Sigma_{y}^{-1} \\\label{eq:optimal-A}
		A &= \Sigma_{x}^{-\half}(\Sigma_x^{\half} (\Sigma_{x} - \Sigma_{xy}\Sigma_{y}^{-1}\Sigma_{xy}^T)\Sigma_{x}^{\half})^\half \Sigma_{x}^{-\half}
	\end{align} 
\end{subequations}
with resulting  transport map
	\begin{equation}\label{eq:T-Gaussian}
		\nabla_x \bar{f}(x,y) = m_x + A(x-m_x) + K(y-m_y).
	\end{equation}
\end{proposition}
\medskip

Note that the above proposition is valid for any joint distribution $P_{XY}$ as long as $\Sigma_X$ is a positive-definite
matrix. However, when $P_{XY}$ is Gaussian, then the map $T$ in \eqref{eq:T-Gaussian} is 
precisely the
linear map that pushes $P_X$ to $P_{Y|X}$. 

The result above can be used to construct a nonlinear filtering algorithm. In particular, when only samples $(X^i_0,Y^i_0) \sim P_{XY}$ are
available, the transport map~\eqref{eq:T-Gaussian} is approximated empirically by substituting the mean and the covariance by their empirical averages, concluding the following update law for the particles: 
\begin{align}\label{eq:OT-EnKF}
X^i_{1}  &= \nabla_x f^{(N)}(X^i_0,y)\nonumber\\&= m_x^{(N)}\! +\! A^{(N)}(X^i_0-m_x^{(N)} ) + K^{(N)}(y- m_y^{(N)})
\end{align}
for any realization of the observation signal $Y=y$, 
where $m_x^{(N)}=N^{-1}\sum_{i=1}^N X^i_0$ and $m_y^{(N)} = N^{-1
} \sum_{i=1}^N Y^i_0$ are the empirical approximation of $m_x$ and $m_y$ respectively, and $K^{(N)}$ and  $A^{(N)}$ are approximations of $K$ and $A$ in~\eqref{eq:optimal-K} and~\eqref{eq:optimal-A} obtained 
from empirical approximations of the covariance matrices.    

The resulting algorithm~\eqref{eq:OT-EnKF} is the optimal transport version of EnKF algorithm for the discrete-time filtering problem. \eqref{eq:OT-EnKF} should be compared with the classical
EnKF update with perturbed observations~\cite{evensen2006,reich11,bergemann2012ensemble}:
	\begin{equation}\label{eq:EnKF}
X^i_1 = X^i_0 + K^{(N)}(y - Y^i_0)
\end{equation}  	   
Interestingly, both \eqref{eq:OT-EnKF} and~\eqref{eq:EnKF} result in the same update rules for the empirical mean and covariance, but different trajectories for the particles. The difference arises due to presence of $A^{(N)}$ in \eqref{eq:OT-EnKF} to ensure that the update law for the particles is optimal with respect to the quadratic  transportation cost. Compared to~\eqref{eq:EnKF},  the optimal transport update~\eqref{eq:OT-EnKF} is expected to be more numerically challenging while admitting lower variance error in the estimation~\cite{taghvaei2020optimal}.  
\begin{remark}
	The update rule~\eqref{eq:OT-EnKF} can also be viewed as the generalization of the continuous-time optimal transport EnKF update of~\cite{taghvaei2020optimal} to the discrete-time setting.
\end{remark}


\subsection{Restriction to ICNNs}\label{sec:ICNN}
The second class of functions discussed here are ICNNs~\cite{AmoXuKol17}. This class of neural networks can be used to represent functions $f(x,y)$ that are convex in $x$.
Universal approximation results have been established for ICNNs 
stating that they can approximate any convex function over a compact domain with a desired accuracy~\cite{CheShiZha18}.
 
 In order to employ ICNNs for the proposed variational problem~\eqref{eq:OT-f}, it is  necessary to represent their convex conjugates.   
Unlike quadratic functions, there are no explicit formulae for the convex conjugates of  ICNNs.
This issue is resolved in~\cite{MakTagOhLee19} by representing the convex conjugate as the 
solution to an inner optimization problem leading to  a min-max problem of the form:
\begin{equation}
\begin{aligned}
&\min_{f \in \text{ICNN}} \max_{g \in \text{ICNN}} \{ \Expect_{P_X \otimes P_Y}[f(X,Y)]  \\&+ \Expect_{P_{XY}}[\nabla_xg(X,Y)^TX - f(\nabla_x g(X,Y),Y)]\} \label{eq:J-ICNN}
\end{aligned}
\end{equation}
The solution to the min-max problem can be numerically approximated using stochastic optimization algorithms resulting in novel nonlinear filtering algorithms for the discrete-time setting. Preliminary numerical results in this direction are presented in Section~\ref{sec:numerics}.

\section{Connection to FPF} \label{sec:FPF}

FPF is an algorithm for the continuous-time nonlinear filtering  problem where the state and observations are modeled as 
continuous-time stochastic processes. For simplicity, assume the state $X\in \Re^n$ is static, while  the observations $\{Z_t \in \Re; ~t \geq 0\}$ is scalar-valued and modeled by the following stochastic differential equation (SDE): 
\begin{equation}\label{eq:obs}
	\ud Z_t = h(X)\ud t + \sigma \ud W_t
\end{equation}
where  $h:\Re^n \to \Re$ is the observation function,  $W_t$ is the standard Wiener process
modeling observation noise, and $\sigma \ge 0$ is the standard deviation of the noise.
In the continuous-time setting, the objective is to compute the conditional distribution $P_{X|\mathcal Z_t}$ where $\mathcal Z_t$
is the filtration generated by the observation process $\{Z_s : ~s\in [0,t]\}$.

The FPF algorithm proceeds by simulating a controlled stochastic process 
\begin{equation}\label{eq:Xbar}
	\ud \bar{X}_t = \k_t(\bar{X}_t) \ud Z_t + u_t(\bar{X}_t) \ud t,\quad \bar{X}_0\sim P_{X}
\end{equation}
where the vector-fields $\k_t$ and $u_t$ are designed  such that the law of $\{\bar{X}_t\}$ coincides with the
posterior distribution $P_{X|\mathcal Z_t}$ for all $t\geq 0$.
	The objective of this section is to formally characterize the vector-fields $\k_t$ and $u_t$ in~\eqref{eq:Xbar}
	directly from the variational formulation~\eqref{eq:OT-f} circumventing application of 
	the continuous-time nonlinear filtering equations for $P_{X|\mathcal Z_t}$ in the original derivation of the FPF algorithm~\cite{yang2016}. For simplicity, the procedure is explained for $t=0$ because the extension to $t>0$ is similar.

	To that end, consider a time discretization of
	the observation process~\eqref{eq:obs} according to
	\begin{equation}\label{eq:obs-FPF}
		Y:= Z_{\Delta t} -Z_{0}= h(X) \Delta t+ \sigma W_{\Delta t}.
	\end{equation}
	The state $X$ and the observation random variable $Y$ are used to define the variational problem~\eqref{eq:OT-f}. 
	 The solution to the variational problem~\eqref{eq:OT-f} is assumed to be of the form
	\begin{equation}\label{eq:f-FPF}
		f(x; y
		) = \frac{1}{2}\|x\|^2 + \phi(x) y + \psi(x)\Delta t,
	\end{equation}
	where $\phi$ and $\psi$ are real-valued functions. The solution can then be used to obtain a transport map  to update $\bar{X}_0$ to $\bar{X}_{\Delta t}$ as follows
\begin{equation}\label{eq:dX}
	\bar{X}_{\Delta t} = \nabla_x f(\bar{X}_0;y) = \bar{X}_0 + \nabla  \phi(\bar{X}_0)y + \nabla \psi(\bar{X}_0)\Delta t
\end{equation} 
which in the limit as $\Delta t\to 0$ yields the stochastic differential equation~\eqref{eq:Xbar} with $\k_0 = \nabla \phi$ and $u_0 = \nabla\psi$.

It remains to identify the functions $\phi$ and $\psi$ by solving the optimization problem~\eqref{eq:OT-f}.  The following proposition identifies the first-order and second-order approximation of the objective function in the asymptotic  limit as $\Delta t \to 0$. 

\begin{proposition}\label{prop:FPF}
	In the asymptotic  limit as $\Delta t \to 0$,  
	the value of the objective function~\eqref{eq:OT-f}, with observation model~\eqref{eq:obs-FPF} and $f$ specified as~\eqref{eq:f-FPF}, is
	\begin{align}\label{eq:J-FPF}
		J(f)&= \frac{1}{2}\sigma_x^2+ J_1(\phi) \Delta t+ J_2(\phi,\psi)\Delta t^2 + O(\Delta t^3) 
	\end{align}
	where $\sigma_x^2 = \Expect [\|X\|^2], 	\hat{h}=\Expect[h(X)]$, and   
		\begin{align*}
			J_1(\phi) &= \Expect\left[\frac{\sigma_w^2}{2}\|\nabla  \phi(X)\|^2 - \phi(X)(h(X)-\hat{h})\right],\\
			J_2(\phi,\psi) &= \Expect\left[\frac{1}{2}\|\nabla \psi(X)-v(X)\|^2 \right] + \tilde{J}_2(\phi), \\ 
				v & = -\frac{h+\hat{h}}{2} \nabla \phi + \frac{\sigma_w^2}{2}\nabla^2 \phi \nabla \phi,
		\end{align*}
	and $\tilde{J}_2$ is only a function of $\phi$ that does not depend on $\psi$. 
\end{proposition}
\medskip

The form of the objective function in~\eqref{eq:J-FPF} suggests that, in the limit as $\Delta t \to 0$,  the minimizer $\phi$ converges to the minimum  of the first-order term $J_1$, and the minimizer $\psi$ converges to the minimum of the second-order term $J_2$ while $\phi$ is fixed. 
Minimizing $J_1$ over $\phi$ yields the first-order optimality condition
\begin{equation}
	-\frac{1}{p(x)} \nabla \cdot (p(x)\nabla \phi(x)) = \frac{1}{\sigma_w^2}(h(x) - \hat{h}),\quad \forall x \in \Re^n,
\end{equation} 
where $p$ is the probability density function of $P_X$. This is known as the weighted Poisson equation. 
Minimizing $J_2$ over $\psi$, while $\phi$ is fixed, yields $\nabla \psi = v + \xi$
where $\xi$ is a divergence-free vector-field , i.e. $\nabla \cdot(p\xi)=0$.

Using the form of $\phi$ and $\psi$ in~\eqref{eq:dX} and taking the continuous-time limit as $\Delta t \to 0$ concludes 
\begin{align*}
	\ud \bar{X}_t =& \nabla \phi_t(\bar X_t) \left(\ud Z_t - \frac{h(\bar{X}_t)+\hat{h}_t}{2}\ud t \right) \\&+  \frac{\sigma_w^2}{2}\nabla^2 \phi_t(\bar X_t) \nabla \phi_t(\bar{X}_t)\ud t + \xi_t(\bar{X}_t)\ud t
\end{align*}   
This is in agreement with the optimal transport form of the FPF algorithm proposed in~\cite{taghvaei2021optimal} and the original form of the  FPF in~\cite{yang2016} modulo the  additional divergence-free term $\xi_t$. 

%
%
\section{Numerical example} \label{sec:numerics}
\begin{figure}[t]
	\centering
	\includegraphics[width=0.9\hsize]{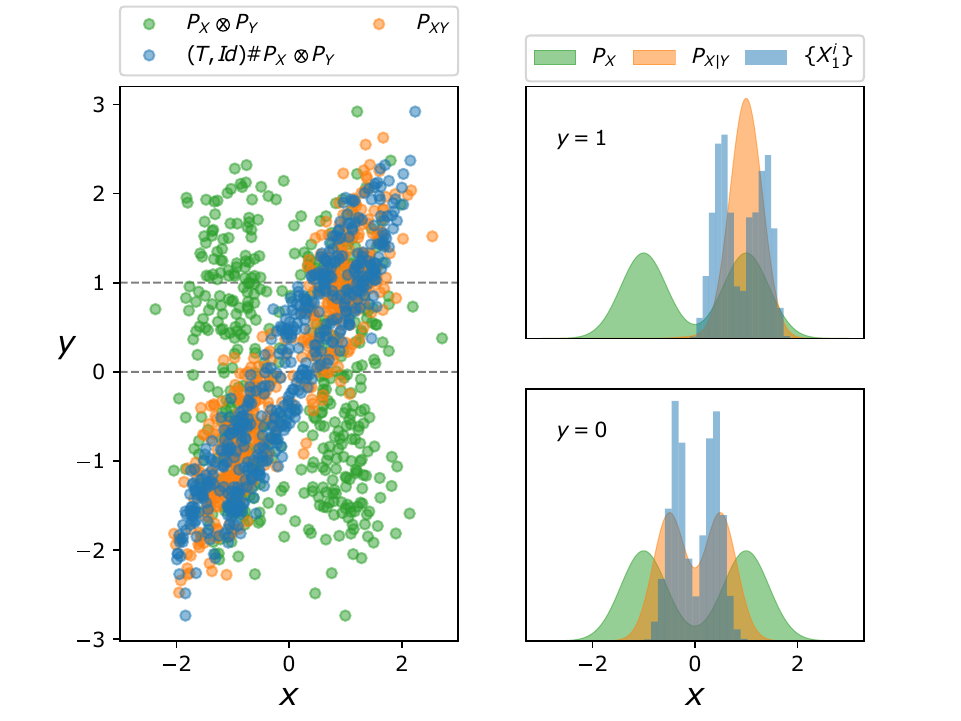}
	\caption{Numerical result for the bimodal example in Section~\ref{sec:numerics}. Samples from the joint distribution $P_{XY}$, the independent coupling $P_X\otimes P_Y$, and the push-forward $(\nabla_x f,\text{Id}) \# P_X\otimes P_Y$ are depicted in the left panel, where $f$ is an approximate solution to the variational problem. problem~\eqref{eq:OT-f}. The prior distribution $P_X$, along with the exact posterior distribution $P_{X|Y=y}$ and the approximated posterior distribution $\nabla_x f(\cdot,y)\#P_X$ for two values $y=0,1$ are depicted in the right panels.   }
	\label{fig:bimodal}
\end{figure}

Let $P_X$ be  a bimodal distribution formed by combination of two Gaussians $N(-1,\sigma^2)$ and $N(+1,\sigma^2)$ with $\sigma^2 = 0.2$.  Assume $Y$ is related to $X$ according to
\begin{equation*}
Y = X + \sigma_w W
\end{equation*} 
where $W$ is standard Gaussian independent of $X$ and $\sigma_w^2 = 0.2$. Consider  solving the variational formulation~\eqref{eq:OT-f} using ICNNs and the min-max formulation~\eqref{eq:J-ICNN}.  For this one dimensional example, consider a simple single layer ICNN architecture where 
\begin{align*}
f(x,y) = \sum_{k=1}^K W_k(W^x_kx+W^y_ky+b_k)_+^2
\end{align*} 
where $W_k\geq 0$, $W^x_k,W^y_k,b_k \in \mathbb R$ for $k=1,\ldots, K$, and $K$ is the size of the network.  A similar architecture is used for $g(x,y)$. Stochastic optimization algorithm (ADAM) is used to solve the min-max problem and learn the parameters of the network. The result is depicted in Figure~\ref{fig:bimodal}. 

\section{Concluding remarks}
The paper presents a variational characterization of the Bayes's law using tools from optimal transportation theory. The variational formulation is used to derive the optimal transport EnKF algorithm, and propose novel generalizations of the EnKF algorithm to the non-Gaussian setting, utilizing ICNN and stochastic optimization algorithms. 
The paper presents preliminary numerical result that serve as proof of concept, while extensive numerical studies and comparison with other nonlinear filtering algorithms are subject of ongoing work.
\label{sec:conclusion}

%

\appendix
\subsection{Proof sketch of the proposition~\ref{prop:quadratic}}
The value of the objective function~\eqref{eq:J-FQ} is obtained  using the quadratic form of the function $f(x,y;\theta)= \frac{1}{2}x^TAx + x^T(Ky+b)$, and its convex conjugate
\begin{align*}
f^*(x,y;\theta) 
=\frac{1}{2}(x-Ky-b)^TA^{-1}(x-Ky-b)
\end{align*} 
The objective function is convex with respect to  $\theta=\{A,K,b\}$ because, (i) the function $f(x,y;\theta)$ is linear in $\theta$, and (ii) $f^*(x,y;\theta)$ is a maximization over linear functions of $\theta$, hence convex.  The solution to the optimization problem is obtained using the first-order optimality condition. 
\subsection{Proof sketch of the proposition~\ref{prop:FPF}}
Express~\eqref{eq:f-FPF} as $f(x,y) =  \frac{1}{2}\|x\|^2  + \eta(x,y)$ where $\eta(x,y)=\phi(x)y + \psi(x)\Delta t$. Then, the convex conjugate of $f$ up to the fourth order in $\eta$ is
\begin{align*}
f^*(x,y) 
&=\frac{1}{2}\|x\|^2 -  \eta(x,y)+ \frac{1}{2}\|\nabla_x \eta(x,y)\|^2 \\&- \frac{1}{2} \nabla_x \eta(x,y)^T \nabla_x^2 \eta(x,y) \nabla_x \eta(x,y)\\
&+ \frac{1}{2} \|\nabla^2_x\eta (x,y) \nabla_x \eta(x,y)\|^2  
\end{align*}
The expression~\eqref{eq:J-FPF}  is obtained using the form for $f$, the expansion of the convex conjugate $f^*$, and the observation model $Y = h(X)\Delta t + \sigma_wW_{\Delta t}$. The details are removed on the account of the space. 

\bibliographystyle{plain}
\bibliography{TAC-OPT-FPF}

\end{document}